\documentclass{amsart}
\usepackage{amssymb,amsfonts,amsmath}

\newtheorem{theorem}{Theorem}

\newtheorem{corollary}[theorem]{Corollary}

\input{tcilatex}

\begin{document}
\title{The formula $ABA=Tr(AB)A$ for matrices}
\author{Grigore C\u{a}lug\u{a}reanu}
\address{Department of Mathematics, Babe\c{s}-Bolyai University, 1 Kog\u{a}%
lniceanu Street, Cluj-Napoca, Romania}
\email{calu@math.ubbcluj.ro}

\begin{abstract}
We prove that this formula characterizes the square matrices over
commutative rings for which all $2\times 2$ minors equal zero.
\end{abstract}

\maketitle

\section{Introduction}

An important notion in Ring Theory is obtained for every idempotent $e$ of a
ring $R$: the set of products $eRe$, called a \textsl{corner }of the ring $R$%
, which is readily seen to be a ring itself with $e$ the multiplicative
identity. The starting point of this note was a Ring Theory research in
which elements of form $eae$ with idempotent $e$ and unit (or nilpotent) $a$
occurred, that is, some specific elements in some corners of a ring. To have
more examples, we tried to check what can be said about matrix products of
this sort.

If one considers $2\times 2$ matrices over commutative rings, it is readily
checked (by direct computation) that if $\det (A)=0$ and $B$ is an arbitrary
matrix then%
\begin{equation*}
ABA=Tr(AB)A
\end{equation*}%
formula we did not find in the literature. A little bit harder, but still
not difficult, the formula can also be checked (by direct computation) for $%
3\times 3$ matrices, whenever all the $2\times 2$ minors of $A$ equal zero.
As for $4\times 4$ matrices (after successfully checking only one entry) we
remembered a good advise given by the 2004 Nobel Prize in Physics, Frank
Wilczek: "I try to avoid hard work. When things look complicate, that is
often a sign that there is a better way to do it".

So we decided to try some different, more conceptual approach in order to
prove the formula for $n\times n$ matrices over commutative rings and any $%
n\geq 2$. Surprisingly, the converse also holds and so we have a
characterization result.

\begin{theorem}
Let $A$, $B$ be $n\times n$ matrices over a commutative ring $R$ and $n\geq 2
$. The formula 
\begin{equation*}
ABA=Tr(AB)A=ATr(BA),
\end{equation*}%
holds for every matrix $B$ if and only if all $2\times 2$ minors of $A$
equal zero.
\end{theorem}

Notice that the second equality follows from the hypothesis that the base
ring is commutative and the well-known trace formula $Tr(AB)=Tr(BA)$.

In section 2 we present the proofs in two special cases: when $A$ has inner
rank one (that is, has a column-row decomposition) and when the entries of $A
$ belong to a field, respectively.

In section 3 we give the proof of the characterization, by induction and
using block multiplication one way, in the general commutative base ring
case. The final section 4 includes some comments and applications.

Recall the following

\textbf{Definition}. The \textsl{inner rank} of an $n\times n$ matrix $A$
over a ring $R$ is the least integer $r$ such that $A$ can be expressed as a
product of an $n\times r$ matrix and an $r\times n$ matrix. For example,
over a division ring this notion coincides with the usual notion of rank. A
square matrix is called \textsl{full} if its inner rank equals its order,
and \textsl{non-full} otherwise.

It is easy to see that, over any commutative ring and for any positive
integer $n\geq 2$, the determinant of an $n$-column-$n$-row product is zero.
Obviously, such products have inner rank 1.

We also recall that a ring was called \textsl{pre-Schreier} (see \cite{zaf})
if every element is \textsl{primal} (i.e., if $r$ divides a product $xy$
then $r=ab$ with $a\mid x$ and $b\mid y$).

\section{Proofs of special cases}

It is easy to prove our formula\textbf{\ if} $A=\mathbf{c}\cdot \mathbf{r}$ 
\textbf{has a column-row decomposition} (i.e., the inner rank of $A$ is 1,
the matrix is \textsl{non-full}) over any commutative ring.

\begin{proof}
A simple proof uses the \textsl{cyclicity of the trace}, i.e., the trace is
invariant with respect to cyclic permutations of its argument. For example, $%
Tr(ABC)=Tr(BCA)=Tr(CAB)$. Let $\mathbf{c}$ be an $n\times 1$ column and let $%
\mathbf{r}$ be an $1\times n$ row such that $A=\mathbf{c}\cdot \mathbf{r}$.
Then $Tr(AB)=Tr(\mathbf{c}\cdot (\mathbf{r}\cdot B))=Tr(\mathbf{r}\cdot
B\cdot \mathbf{c})=\mathbf{r}\cdot B\cdot \mathbf{c}$ since this is a
scalar. Therefore $ABA=\mathbf{c}\cdot \mathbf{r}\cdot B\cdot \mathbf{c}%
\cdot \mathbf{r}=\mathbf{c}Tr(AB)\mathbf{r}=Tr(AB)A$.
\end{proof}

However, having a column-row decomposition for any matrix with all $2\times
2 $ minors equal to zero, requires additional conditions on the base ring.

For instance, in the $n=2$ case, it is proved (see \cite{CP}) that the
existence of a column-row decomposition for a zero determinant $2\times 2$
matrix is equivalent to the base ring being a pre-Schreier domain.

As it is well-known, the proof is valid if the base ring is a field, as, for
any matrix, the inner rank equals the rank.

\bigskip

Next we supply \textbf{a different proof for }(square) \textbf{matrices over
fields}.

First we mention two facts which are easier to describe for linear maps
(instead of matrices).

Suppose $f$ is a linear map from $V$ to itself, where $V$ is some finite $n$%
-dimensional vector space over a field $K$. Recall that, by definition, the
rank of $f$ is $r=\dim (f(V))$.

\textbf{1}. Suppose $r=1$. Then the range $f(V)$ has a vector basis, say $\{%
\mathbf{v}\}$ for some $\mathbf{v}\neq \mathbf{0}$. Since for every $\mathbf{%
x}\in V$, $f(\mathbf{x})\in f(V)$ it follows that $f(\mathbf{x})=d\mathbf{v}$
for some $d\in K$. In particular $f(\mathbf{v})=c\mathbf{v}$ for some $c\in K
$, so $\mathbf{v}$ is an eigenvector corresponding to the eigenvalue $c$.

\textbf{2}. Suppose $r=1$ and so the nullity $\dim (\ker (f))=n-1$. Since
the multiplicity of an eigenvalue is at least the dimension of the
corresponding eigenspace, it follows that $0$ is an eigenvalue with
multiplicity at least $n-1$. Moreover, as the sum of all eigenvalues
(counted with multiplicity) is $Tr(f)$, the last eigenvalue is $Tr(f)$.

Now, given some square matrix $A$, we can apply this to the map canonically
attached to $A$. Since the formula trivially holds for $B=0$, we assume $%
B\neq 0$.

\begin{proof}
Suppose $rank(A)=1$ (and so $A\neq 0$). By the above (see \textbf{1}), there
is a nonzero vector $\mathbf{v}$ such that for each vector $\mathbf{x}$, $A%
\mathbf{x}=d\mathbf{v}$ for some $d\in R$. Then we compute%
\begin{equation*}
ABA\mathbf{x}=AB(d\mathbf{v})=dA(B\mathbf{v})=dc\mathbf{v}=c(d\mathbf{v})=cA%
\mathbf{x}
\end{equation*}%
for some $c\in R$ (that is, $A(B\mathbf{v})=c\mathbf{v}$). Since the
equality holds for every $\mathbf{x}$, $ABA=cA$ follows.

Finally, it remains to notice that $rank(AB)=1$ and since $(AB)\mathbf{v}=c%
\mathbf{v}$ we deduce (see \textbf{2} above) that $c=Tr(AB)$.
\end{proof}

Despite the fact that notions like rank of a matrix, linearly independent
vector (or column in $R^{n}$), null space, eigenvalue and eigenvector can be
defined for matrices over commutative rings (see \cite{bro}), this proof
cannot be extended to the general case of matrices over commutative rings.
There are several impediments occurring.

One of these is that if $c\in R$ is an eigenvalue (e.g., a root of the
characteristic polynomial $p_{A}(X)$) it does not follow that the
eigenvector of $A$ associated to $c$ is independent over $R$. Of course, $%
p_{A}(X)$ may not have roots in $R$.

Moreover, $rank(A)=1$ implies that the $2\times 2$ minors are zero, but not
conversely. Unlike the field case, a matrix can have rank zero without being
the zero matrix.

\section{Proof of the theorem}

The formula is obviously true for $n=1$ without any hypothesis. Indeed $%
aba=(ab)a=a(ba)$ holds over any (associative) ring (possibly not commutative
nor unital).

\begin{proof}
One way, using induction and block multiplication, we provide a proof of the
formula for $n\times n$ matrices over commutative rings, whenever all $%
2\times 2$ minors equal zero. As noticed in the Introduction, if $\det (A)=0$
for a $2\times 2$ matrix (over any commutative ring) $A=[a_{ij}]$ and $%
B=[b_{ij}]$ is arbitrary, it is not hard to check $ABA=Tr(AB)A$ (here $%
Tr(AB)=a_{11}b_{11}+a_{12}b_{21}+a_{21}b_{12}+a_{22}b_{22}$). 

Let $A,B\in \mathbb{M}_{n}(R)$ for a commutative ring $R$. In both $A$ and $B
$ we emphasize the $(n-1)\times (n-1)$ left upper corner, that is, we write $%
A=\left[ 
\begin{array}{cc}
A^{\prime } & \mathbf{a}_{2} \\ 
\mathbf{a}_{1} & a_{nn}%
\end{array}%
\right] $ with $A^{\prime }\in \mathbb{M}_{n-1}(R)$, one $1\times (n-1)$ row 
$\mathbf{a}_{1}=\left[ 
\begin{array}{ccc}
a_{n1} & ... & a_{n,n-1}%
\end{array}%
\right] $ and one $(n-1)\times 1$ column $\mathbf{a}_{2}=\left[ 
\begin{array}{c}
a_{1n} \\ 
\vdots  \\ 
a_{n-1,n}%
\end{array}%
\right] $. We use similar notations for $B$ and block multiplication for
these matrices.

Then $AB=\left[ 
\begin{array}{cc}
A^{\prime }B^{\prime }+\mathbf{a}_{2}\mathbf{b}_{1} & A^{\prime }\mathbf{b}%
_{2}+\mathbf{a}_{2}b_{nn} \\ 
\mathbf{a}_{1}B^{\prime }+a_{nn}\mathbf{b}_{1} & \mathbf{a}_{1}\mathbf{b}%
_{2}+a_{nn}b_{nn}%
\end{array}%
\right] $ and $ABA=$

\noindent \noindent $\left[ 
\begin{array}{cc}
(A^{\prime }B^{\prime }+\mathbf{a}_{2}\mathbf{b}_{1})A^{\prime }+(A^{\prime }%
\mathbf{b}_{2}+\mathbf{a}_{2}b_{nn})\mathbf{a}_{1} & (A^{\prime }B^{\prime }+%
\mathbf{a}_{2}\mathbf{b}_{1})\mathbf{a}_{2}+(A^{\prime }\mathbf{b}_{2}+%
\mathbf{a}_{2}b_{nn})a_{nn} \\ 
(\mathbf{a}_{1}B^{\prime }+a_{nn}\mathbf{b}_{1})A^{\prime }+(\mathbf{a}_{1}%
\mathbf{b}_{2}+a_{nn}b_{nn})\mathbf{a}_{1} & (\mathbf{a}_{1}B^{\prime
}+a_{nn}\mathbf{b}_{1})\mathbf{a}_{2}+(\mathbf{a}_{1}\mathbf{b}%
_{2}+a_{nn}b_{nn})a_{nn}%
\end{array}%
\right] $. By induction hypothesis, assume $A^{\prime }B^{\prime }A^{\prime
}=Tr(A^{\prime }B^{\prime })A^{\prime }$. In order to prove $%
ABA=Tr(AB)A=Tr(AB)\left[ 
\begin{array}{cc}
A^{\prime } & \mathbf{a}_{2} \\ 
\mathbf{a}_{1} & a_{nn}%
\end{array}%
\right] $ we have to check four equalities:

(i) $Tr(A^{\prime }B^{\prime })A^{\prime }+\mathbf{a}_{2}\mathbf{b}%
_{1}A^{\prime }+(A^{\prime }\mathbf{b}_{2}+\mathbf{a}_{2}b_{nn})\mathbf{a}%
_{1}=Tr(AB)A^{\prime }$

(ii) $(A^{\prime }B^{\prime }+\mathbf{a}_{2}\mathbf{b}_{1})\mathbf{a}%
_{2}+(A^{\prime }\mathbf{b}_{2}+\mathbf{a}_{2}b_{nn})a_{nn}=Tr(AB)\mathbf{a}%
_{2}$

(iii) $(\mathbf{a}_{1}B^{\prime }+a_{nn}\mathbf{b}_{1})A^{\prime }+(\mathbf{a%
}_{1}\mathbf{b}_{2}+a_{nn}b_{nn})\mathbf{a}_{1}=Tr(AB)\mathbf{a}_{1}$

(iv) $(\mathbf{a}_{1}B^{\prime }+a_{nn}\mathbf{b}_{1})\mathbf{a}_{2}+(%
\mathbf{a}_{1}\mathbf{b}_{2}+a_{nn}b_{nn})a_{nn}=Tr(AB)a_{nn}$.

Since $Tr(AB)=Tr(A^{\prime }B^{\prime }+\mathbf{a}_{2}\mathbf{b}_{1})+%
\mathbf{a}_{1}\mathbf{b}_{2}+a_{nn}b_{nn}$, and $Tr(\mathbf{a}_{2}\mathbf{b}%
_{1})=Tr(\mathbf{b}_{1}\mathbf{a}_{2})=\mathbf{b}_{1}\mathbf{a}_{2}$, the
equalities amount to

(i) $\mathbf{a}_{2}\mathbf{b}_{1}A^{\prime }+(A^{\prime }\mathbf{b}_{2}+%
\mathbf{a}_{2}b_{nn})\mathbf{a}_{1}=(\mathbf{b}_{1}\mathbf{a}_{2}+\mathbf{a}%
_{1}\mathbf{b}_{2}+a_{nn}b_{nn})A^{\prime }$

(ii) $(A^{\prime }B^{\prime }+\mathbf{a}_{2}\mathbf{b}_{1})\mathbf{a}%
_{2}+A^{\prime }\mathbf{b}_{2}a_{nn}=(Tr(A^{\prime }B^{\prime })+\mathbf{b}%
_{1}\mathbf{a}_{2}+\mathbf{a}_{2}\mathbf{b}_{1})\mathbf{a}_{2}$

(iii) $(\mathbf{a}_{1}B^{\prime }+a_{nn}\mathbf{b}_{1})A^{\prime
}=(Tr(A^{\prime }B^{\prime })+\mathbf{b}_{1}\mathbf{a}_{2})\mathbf{a}_{1}$

(iv) $\mathbf{a}_{1}B^{\prime }\mathbf{a}_{2}=Tr(A^{\prime }B^{\prime
})a_{nn}$.

We just provide some details of how all these equalities are verified.

(i) The LHS has a sum of three $(n-1)\times (n-1)$ matrices and the RHS has $%
A^{\prime }$ multiplied by a scalar. We have to check the equalities of the
entries in the LHS matrix respectively in the RHS matrix. Each entry is a
linear combination of $b$'s with coefficients products of two $a$'s. Some
corresponding entries are already equal (e.g., coefficients of $b_{n1}$ in
the corresponding $(1,n-1)$ entries), the other use the vanishing of the $%
2\times 2$ minors (e.g., we use $\left\vert 
\begin{array}{cc}
a_{n-1,1} & a_{n-1,n-1} \\ 
a_{n1} & a_{n,n-1}%
\end{array}%
\right\vert =0$ in order to check the equality of the products which
multiply $b_{n-1,n}$, for the $(n-1,1)$ entries).

(ii) The LHS is a sum of three columns and the RHS is a product of a scalar
and a column. To check the equality amounts to verify the equality of the
corresponding $n-1$ entries. As in the previous case, each entry is a linear
combination of $b$'s with coefficients products of two $a$'s. Some
corresponding entries are already equal (e.g., coefficients of $b_{n1}$ in
the corresponding upper entries), the other use the vanishing of the $%
2\times 2$ minors (e.g., we use $\left\vert 
\begin{array}{cc}
a_{11} & a_{1n} \\ 
a_{n-1,1} & a_{n-1,n}%
\end{array}%
\right\vert =0$ in order to check the equality of the products which
multiply $b_{1,n-1}$, for the upper entries).

(iii) Both the LHS and RHS of the equality are linear combination of all $%
b_{ij}$ for $1\leq i\leq n,1\leq j\leq n-1$ with coefficients products of
two $a$'s. For every pair $(i,j)$, the equality of the coefficients of $%
b_{ij}$ in both sides amounts to the vanishing of a minor of type $%
\left\vert 
\begin{array}{cc}
a_{ji} & a_{jj} \\ 
a_{ni} & a_{nj}%
\end{array}%
\right\vert $ or $\left\vert 
\begin{array}{cc}
a_{j1} & a_{ji} \\ 
a_{n1} & a_{ni}%
\end{array}%
\right\vert $ or similar.

(iv) Both the LHS and RHS of the equality are linear combination of all $%
b_{ij}$ for $1\leq i,j\leq n-1$ with coefficients products of two $a$'s. For
every pair $1\leq i,j\leq n-1$, the equality of the coefficients of $b_{ij}$
in both sides amounts precisely to the vanishing of the minor $\left\vert 
\begin{array}{cc}
a_{ji} & a_{jn} \\ 
a_{ni} & a_{nn}%
\end{array}%
\right\vert $.

As for the converse, take an arbitrary minor $m=\left\vert 
\begin{array}{cc}
a_{ik} & a_{il} \\ 
a_{jk} & a_{jl}%
\end{array}%
\right\vert $ and choose $B=E_{lj}$. It is easy to see that $C=:AE_{lj}A=%
\mathrm{col}_{l}(A)\cdot \mathrm{row}_{j}(A)$ and since $Tr(AE_{lj})=a_{jl}$%
, $D=:Tr(AE_{lj})A=a_{jl}A$. Now compute $c_{ik}=a_{jk}a_{il}$ and $%
d_{ik}=a_{jl}a_{ik}$. Since $C=D$, the minor $m$ indeed equals zero.
\end{proof}

\textbf{Remark}. For $n=2$, the converse (that is, $\det (A)=0$) follows
easily by the Cayley-Hamilton theorem, taking $B=I_{2}$. In the $n=3$ case,
the same choice, $B=I_{3}$, and the Cayley-Hamilton theorem still gives $%
\det (A)=0$ (since $A^{2}=Tr(A)A$ and $Tr(A^{2})=Tr^{2}(A)$) but not the
vanishing of all the $2\times 2$ minors.

\section{Comments and Applications}

Two special cases of this formula can be found on MathOverflow.

1) Let $A$ be an $n\times n$ complex matrix having rank 1. Prove that $%
A^{2}=cA$ for some scalar $c$ (see \cite{lee}, solution by M. van Leeuwen
2014), and

2) Show that if $Tr(AB)=0$ and $A$ has rank 1 then $ABA=0$ (see \cite{euy},
solution by EuYu 2016).

The first follows from our formula by taking $B=I_{n}$, where $c$ turns out
to be precisely $Tr(AB)$, and the second follows directly from the formula,
whenever $Tr(AB)=0$.

However, as our theorem proves, both hold over any commutative ring, not
only for complex matrices, or over special integral domains (where a
column-row splitting is possible).

An easy but general case which was not mentioned follows from
Cayley-Hamilton theorem: let $A$ be a zero determinant $2\times 2$ matrix
over any commutative ring. Then $A^{2}=Tr(A)A$.

\bigskip

The following applications are straightforward.

\begin{corollary}
Assume $n\geq 2$, all $2\times 2$ minors of the $n\times n$ matrix $A$ equal
zero and the $n\times n$ matrix $B$ is arbitrary. Then

(a) $(AB)^{2}=Tr(AB)AB$ but not conversely, unless $B$ is a unit.

(b) $Tr(ABA)=Tr(AB)Tr(A)=Tr(A)Tr(BA)$; in particular $Tr(A^{2})=Tr^{2}(A)$.
\end{corollary}

It is well-known that if $R$ is any commutative ring, the endomorphisms of
an $R$-module $M$ form an algebra over $R$ denoted $End_{R}(M)$. Then there
is a canonical $R$-linear map: $M^{\ast }\otimes _{R}M\longrightarrow R$
induced through linearity by $f\otimes x\mapsto f(x)$; it is the unique $R$%
-linear map corresponding to the natural pairing. If $M$ is a finitely
generated projective $R$-module, then one can identify $M^{\ast }\otimes
_{R}M=End_{R}(M)$ through the canonical homomorphism mentioned above and
then the above is the trace map: $Tr:End_{R}(M)\longrightarrow R$. Our
formula can be transferred \emph{mutatis mutandis} to this context.

\end{document}